\documentclass[12pt]{article}

\usepackage[T1]{fontenc}
\usepackage[latin1]{inputenc}
\usepackage[british]{babel}
\usepackage{amsfonts}
\usepackage{amsthm}
\usepackage{amsmath}

\DeclareTextCompositeCommand{\v}{OT1}{t}{%
    t\kern-.23em\raise.24ex\hbox{'}}

\newtheorem{proposition}{Proposition}
\newtheorem{lemma}[proposition]{Lemma}

\newcommand{\C}{\ensuremath{\mathbb C}}
\newcommand{\R}{\ensuremath{\mathbb R}}
\newcommand{\Z}{\ensuremath{\mathbb Z}}
\newcommand{\N}{\ensuremath{\mathbb N}}

\newcommand{\HH}{\ensuremath{{\mathcal H}}}
\newcommand{\UU}{\ensuremath{{\mathcal U}}}

\newcommand{\gs}{\ensuremath{{\mathfrak s}}}
\newcommand{\gl}{\ensuremath{{\mathfrak l}}}
\newcommand{\go}{\ensuremath{{\mathfrak o}}}

\newcommand{\Uqsl}{\UU_{q}(\gs\gl(2,\R))}

\DeclareMathOperator{\mspan}{span}
\DeclareMathOperator{\LHS}{LHS}
\DeclareMathOperator{\RHS}{RHS}

\renewcommand{\imath}{\mathrm{i}}

\begin{document}

\title{Construction of raising and lowering operators for $\Uqsl$}

\author{Pavel \v{S}\v{T}OV\'I\v{C}EK\\
  \\
  Department of Mathematics\\
  Faculty of Nuclear Science, Czech Technical University\\
  Trojanova 13, 120 00 Prague, Czech Republic\\
  E-mail: stovicek@KMLinux.fjfi.cvut.cz}
\date{}


\maketitle

\abstract{\noindent
  In the case of $\UU_q(\gs\gl(2,\R))$ at root of unity
  $q$-deformed analogues are proposed for the generator of the maximal
  compact subalgebra, $J$, and for the raising and lowering operators.
  We prove an algebraic identity which implies that $J$ has similar
  properties as in the nondeformed case.}

\section{Introduction}
\label{sec:Introduction}

Despite of an intensive development of the theory of quantum groups
during the last decade the results concerning real noncompact forms
are rather rare and far away of being complete. Even in the lowest
rank cases one encounters serious difficulties (see, for example,
references \cite{[1],[2]}). Naturally, most attention has been paid to
groups with a physical interpretation \cite{[3],[11]}. Here we
concentrate on the quantum real group $\UU_q(\gs\gl(2,\R))$. Even this
simple example lacks, to the author's knowledge, a consistent theory
on the $q$-deformed level.

Let us first recall some basic facts concerning the irreducible
unitary representations of the classical Lie group $SL(2,\R)$
\cite{[4],[5]}. Let $H, E, F$ be the standard basis elements of
$\gs\gl(2,\R)$. All of them are skew self-adjoint with respect to the
corresponding $\ast$-involution in $\UU(\gs\gl(2,\C))$. Set
\begin{equation*}
  J=E-F,\ W^\pm = H\pm\imath(E+F).
\end{equation*}
For any irreducible representation $\varrho$ of $SL(2,\R)$ in a
Hilbert space $\HH$ the spectrum of $-\imath\varrho(J)$ is simple and
contained in $\Z$. If $\{v_n\}_{n\in M}$, $M\subset\Z$, is a basis of
$\HH$ formed by suitably normalised eigenvectors of
$-\imath\varrho(J)$ then
\begin{equation*}
  -\imath\varrho(J)v_n=n\,v_n,\
  \varrho(W^\pm)v_n=(s+1\pm n)v_{n\pm2},
\end{equation*}
where $s\in\C$ is a fixed parameter.

Our main goal here is to propose $q$-deformed analogues for the
generator of the compact subalgebra $\gs\go(2)$, $J$, and for the
raising and lowering operators $W^\pm$. This short contribution does
not aim at all, however, to develop a complete representation theory
for $\UU_q(\gs\gl(2,\R))$; this task is deferred to another
publication.

\section{$\UU_q(\gs\gl(2,\R))$ at root of unity}
\label{sec:Uq(sl(2,R))}

We shall consider the case when $q$ is a root of unity. We assume that
\begin{equation*}
  q=\exp(2\pi\imath  P/Q)
\end{equation*}
where $Q\in \N $ is odd, $P\in \{1,2,\ldots,Q-1\}$ and $P$ and $Q$ are
relatively prime integers. Let $K, K^{-1}, E, F$ be the standard
generators of the complex Hopf algebra $\UU_q(\gs\gl(2,\C))$. The
$\ast$-involutions defining the real form is prescribed by the
relations
\begin{equation*}
  K^{\ast }=K,\textrm{ }E^{\ast }=-q^{-1}E,\textrm{ }F^{\ast }=-q\, F.
\end{equation*}

Here we prefer another set of generators generating a ``slightly
smaller'' Hopf algebra $\UU$, namely
\begin{equation*}
  X=-\imath\, q^{-1}E\, K^{-1},\textrm{ }Y=-\imath\, q\, F\, K^{-1},
  \textrm{ }Z=K^{-2}.
\end{equation*}
Thus the defining relations for $\UU$ are
\begin{equation}
  \label{eq:defrel}
  \begin{aligned}
    {}& Z {Z^{-1}}={Z^{-1}}Z=1,\ Z X={q^{-2}}X Z,\ Z\, Y={q^2}\,Y Z,  \\
    {}& {q^{-1}}X Y-q\, Y X= \frac{1}{q-{q^{-1}}} ({Z^2}-1),
  \end{aligned}
\end{equation}
the comultiplication is defined by
\begin{equation*}
  \Delta Z=Z\otimes Z,\ \Delta X=1\otimes X+X\otimes Z,\
  \Delta Y=1\otimes Y+Y\otimes Z,
\end{equation*}
the antipode is defined by
\begin{equation*}
  S(Z)={Z^{-1}},\ S(X)=-X {Z^{-1}},\ S(Y)=-Y {Z^{-1}},
\end{equation*}
and the counit is defined by
\begin{equation*}
  \varepsilon(Z)=1,\ \varepsilon(X)=\varepsilon (Y)=0\,.
\end{equation*}
The $\ast$-involution defining the real form now reads
\begin{equation}
  \label{eq:starorig}
  Z^\ast=Z,\ X^\ast=X,\ Y^\ast=Y\,.
\end{equation}


The irreducible representations are divided into two families (see the
original paper \cite{[6]} or the monographes \cite{[7],[8]}).
The representations from the first family are labelled by a discrete
index $r=0,1,\ldots,Q-1$, the representation space is
$\HH_r=\mspan\{{v_0},{v_1},\ldots,{v_r}\}$, and the representation
is defined by
\begin{equation}
  \label{eq:family1}
  Z\cdot{v_j}=\pm {q^{r-2j}}{v_j},\
  X\cdot{v_j}=-{{{q^{r-2j-1}}[r-j]}_q}{v_{j+1}},\
  Y\cdot{v_j}={{[j]}_q}{v_{j-1}}.
\end{equation}
The representations from the second family are labelled by 3 complex
parameters $\lambda, a, b$, the representation space is always
$\HH=\{{v_0},{v_1},\ldots,{v_{Q-1}}\}$, the representation is defined
by
\begin{equation}
  \label{eq:family2}
  \begin{aligned}
    {}& Z\cdot{v_j}=\lambda\, {q^{2j}}{v_j},  \\
    {}& X\cdot{v_j}=-\imath\, {q^{j-1}}\bigg(a b-{[j]_q}
    \frac{\lambda\, q^{j-1} - \lambda^{-1} q^{-j+1}}{q-{q^{-1}}}\bigg)
    {v_{j-1}}\,\quad \textrm{for } j\neq 0,  \\
    {}& X\cdot{v_0}=-\imath\, {q^{-1}}a\, {v_{Q-1}},  \\
    {}& Y\cdot{v_j}=-\imath\, \lambda\, {q^{j+1}}{v_{j+1}}
    \quad\textrm{for } j\neq Q-1,\\
    {}& Y\cdot{v_{Q-1}}=-\imath\, \lambda\, b\, {v_0}.
  \end{aligned}
\end{equation}
The two families intersect for the values of parameters $r=Q-1$ in
(\ref{eq:family1}) and $a=b=0$, $\lambda =\pm q^{-Q+1}$ in
(\ref{eq:family2}).

Unfortunately, as pointed out in reference \cite{[9]}, the only
representations from this list which are unitarizable with respect to
the $\ast$-involution (\ref{eq:starorig}) are one-dimensional. This
absolutely unsatisfactory situation of course appeals for a modified
definition of $\UU_q(\gs\gl(2,\R))$. A~possible way out has been
proposed in reference \cite{[9]}. The set of generators is augmented
by introducing an additional generator $T$ which satisfies
\begin{equation*}
  T^{2}=1,\textrm{ }\Delta T=T\otimes T,\textrm{ }S(T)=T,
  \textrm{ }\varepsilon(T)=1.
\end{equation*}
The augmented Hopf algebra will be called $\UU^\natural$. The modified
$\ast$-involution takes the form
\begin{equation}
  \label{eq:starmod}
  X^{\ast }=T\, X\, T,\textrm{ }Y^{\ast }=
  T\, Y\, T,\textrm{ }Z^{\ast }=T\, Z\, T,\textrm{ }T^{\ast }=T.
\end{equation}

As shown in papers \cite{[9],[10]}, all representations from the first family
(\ref{eq:family1}) are unitarizable in this modified sense and may be
regarded as $q$-deformations of representations of $SL(2,\R)$
belonging to the discrete series.

\section{Raising and lowering operators}

Set
\begin{equation*}
  J=(q X-{q^{-1}}Y){Z^{-1}}={Z^{-1}}({q^{-1}} X-q Y) .
\end{equation*}
One easily finds that
\begin{equation}
  \label{eq:Jrelations}
  \Delta J={Z^{-1}}\otimes J+J\otimes 1,\ S(J)=-Z\,J,\
  \varepsilon(J)=1\, .
\end{equation}
For the extended algebra $\UU^\natural$ it seems to be reasonable to
impose the condition
\begin{equation*}
  J T=T J,
\end{equation*}
i.e. ${J^\ast}=J$.

Note that conversely one can express
\begin{equation*}
  X=\frac{1}{{q^2}-{q^{-2}}}(q J Z-{q^{-1}}Z J),\
  Y=\frac{1}{{q^2}-{q^{-2}}}({q^{-1}} J Z-q Z J)\,.
\end{equation*}
It is also straightforward to derive the relations between $J$ and
$Z$. The second and third relations in (\ref{eq:defrel}) give
\begin{equation}
  \label{eq:ZJrelation1}
  {Z^2}J-({q^2}+{q^{-2}})Z J Z+J {Z^2}=0
\end{equation}
while the last one implies
\begin{equation}
  \label{eq:ZJrelation2}
  ({q^2}+1+{q^{-2}})\, Z {J^2}Z-J Z J Z-J {Z^2}J-Z J Z J
  ={{(q+{q^{-1}})}^2}({Z^2}-1)\,.
\end{equation}
All the above manipulations are reversible and this shows
\begin{proposition}
  The Hopf algebra $\UU$ is isomorphic to the Hopf algebra $\UU'$
  generated by $Z$ and $J$ with the defining relations
  (\ref{eq:ZJrelation1}), (\ref{eq:ZJrelation2}),
  (\ref{eq:Jrelations}) and the known relations for $Z$
  ($Z\,Z^{-1}=Z^{-1}Z=1$, $\Delta Z=Z\otimes Z$, $S(Z)=Z^{-1}$,
  $\varepsilon(Z)=1$).
\end{proposition}

The quantum numbers $[x]_q=(q^x-q^{-x})/(q-q^{-1})$ are supposed to
be defined for any $x\in\C$ for we set, by definition,
$q^x=\exp(2\pi\imath(P/Q)x)$. The main result of this contribution is
the following identity.

\begin{proposition}
  For all $x\in\C$ it holds true that
  \begin{equation}
    \label{eq:identity}
    \begin{aligned}
      {}& Z\big(J-{[x+2]_q}\big)\big(J-{[x]_q}\big)
      \big(J-{[x-2]_q}\big)Z \\
      {}& \quad
      =\Big(\big(J-{[x]_q}\big)Z\big(J-{[x]_q}\big)Z
      -{{{[2]_q}}^2}\Big)\big(J-{[x]_q}\big)\,,
    \end{aligned}
  \end{equation}
\end{proposition}

\begin{proof}
  To prove (\ref{eq:identity}) we expend the both sides in the
  variable $y={q^x}$,
  \begin{equation*}
    \LHS -\RHS =c_{-2}\,y^{-2}+c_{-1}\,y^{-1}+c_0+c_1\, y+c_2\, y^2,
  \end{equation*}
  and we find that
  \begin{equation*}
    \begin{aligned}
      {}& {{(q-{q^{-1}})}^2}c_{-2}={{(q-{q^{-1}})}^2}c_2
      =-J {Z^2}+({q^2}+{q^{-2}}) Z J Z-{Z^2}J,  \\
      {}& (q-{q^{-1}})c_{-1}=-(q-{q^{-1}})c_1 \\
      {}& \qquad\quad
      =({q^2}+1+{q^{-2}})Z {J^2}Z-J Z J Z-J {Z^2}J-Z\, J Z
      J-{{{{[2]}_q}}^2}({Z^2}-1),
    \end{aligned}
  \end{equation*}
  and
  \begin{equation*}
    \begin{aligned}
      {{(q-{q^{-1}})}^2}c_0 =&
      \,2\, (J {Z^2}-({q^2}+{q^{-2}}) Z J Z+{Z^2}J) \\
      & +{{(q-{q^{-1}})}^2}\big( Z {J^3}Z-{{{{[2]}_q}}^2}Z\, J Z-J Z J Z
      J+{{{{[2]}_q}}^2}J\big) .
    \end{aligned}
  \end{equation*}
  Thus the coefficients $c_{\pm2}$, $c_{\pm1}$ vanish owing to
  (\ref{eq:ZJrelation1}) and (\ref{eq:ZJrelation2}), the coefficient
  $c_0$ vanishes if it is true in the particular case $x=0$. The
  following lemma concludes the proof.
\end{proof}

\begin{lemma}
  It holds true that
  \begin{equation*}
    Z\big({J^2}-{{(q+{q^{-1}})}^2}\big)J Z
    =\big(J Z J Z-{{(q+{q^{-1}})}^2} \big)J\,.
  \end{equation*}
\end{lemma}

\begin{proof}
  To prove the particular case of (\ref{eq:identity}) with $x=0$ set
  \begin{equation*}
    V=Z {J^3} Z-{{{{[2]}_q}}^2}Z J Z+{{{{[2]}_q}}^2} J-J Z J Z\,J.
  \end{equation*}
  Then it holds
  \begin{equation*}
    \begin{aligned}
      {}& Z V+V Z = \\
      {}& \quad\big(({q^2}+1+{q^{-2}})\,
      Z J^2 Z-J Z^2 J-Z J Z J-J Z J Z+
      {{(q+{q^{-1}})}^2}(1-{Z^2})\big) \\
      {}& \quad\times J Z + Z J  \\
      {}& \quad\times \big(({q^2}+1+{q^{-2}})\, Z J^2 -J {Z^2} J-Z J Z
      J-J Z J Z
      +{{(q+{q^{-1}})}^2}(1-{Z^2})\big)  \\
      {}& \quad +(J {Z^2} -({q^2}+{q^{-2}}) Z J Z+{Z^2} J) {J^2} Z \\
      {}& \quad +Z {J^2}(J Z^2-({q^2}+{q^{-2}})\, Z J Z+{Z^2} J)
    \end{aligned}
  \end{equation*}
  and thus $Z\, V+V Z=0$. This clearly implies
  \begin{equation*}
    {Z^k}V={{(-1)}^k}\,V Z, \quad\textrm{for all } k\in \Z .
  \end{equation*}
  In particular, ${Z^Q}V=-V {Z^Q}$. But ${Z^Q}$ belongs to the centre
  of the algebra $\UU$ and therefore ${Z^Q}V=V {Z^Q}$. Hence
  ${Z^Q}V=0$ and consequently $V=0$.
\end{proof}

Some straightforward consequences follow from identity
(\ref{eq:identity}). Consider an irreducible representation of $\UU$
whose dimension equals $d$. If $v$ is an eigenvector of $J$ with an
eigenvalue $[x]_q$ then the vectors
\begin{equation}
  \label{eq:spectrum}
  \big(J-{[x]_q}\big)\big(J-{[x-2]_q}\big)Z\cdot v,\
  \big(J-{[x]_q}\big)\big(J-{[x+2]_q}\big)Z\cdot v,
\end{equation}
either vanish or are eigenvectors of $J$ corresponding to the
eigenvalues $[x+2]_q$ and $[x-2]_q$, respectively. The spectrum of $J$
is of the form
\begin{equation*}
  [x]_q, [x+2]_q,\ldots, [x+2d-2]_q,
\end{equation*}
and hence the operator $J$ has similar properties (up to the factor
$\imath$) as its nondeformed counterpart. For example, in the case of
the representations from the first family (\ref{eq:family1}) the
operator $J$ is known to have the spectrum (\ref{eq:spectrum}) with
$x=Q-d+1$ \cite{[9]}. Furthermore, the operators
$(J-{[x]_q})(J-{[x-2]_q})Z$ and $(J-{[x]_q})(J-{[x+2]_q})Z$ play the
role of raising and lowering operators, respectively. Finally, the
matrix $Z$, if expressed in the corresponding eigenbasis of $J$, is
tridiagonal in the case of the first family of representations
(\ref{eq:family1}) or tridiagonal with the cyclic convention modulo
$Q$ in the case of the second family (\ref{eq:family2}).

\section*{Acknowledgements}
The research was partially supported by Grant GACR 201/01/01308.


\begin{thebibliography}{99}

\bibitem{[1]} S.L. Woronowicz: Lett. Math. Phys. {\bf 23} (1991) 251.

\bibitem{[2]} L.I. Korogodsky: Commun. Math. Phys. {\bf 163} (1994) 433.

\bibitem{[3]} V.K. Dobrev, P. Moylan: Phys. Lett. B {\bf 315} (1993) 292.

\bibitem{[11]} V.K. Dobrev and R. Floreanini: J. Phys. {\bf A27} (1994) 4831.

\bibitem{[4]} V. Bargmann: Ann. Math. {\bf48} (1947) 568.

\bibitem{[5]} S. Lang: {\it $SL_2(\mathbf{R})$}.
  Addison-Wesley, Reading, 1975.

\bibitem{[6]} P. Roche, D. Arnaudon: Lett. Math. Phys. {\bf17} (1989) 295.

\bibitem{[7]} V. Chari, A. Pressley:
  {\it A Guide to Quantum Groups}.
  Cambridge University Press, Cambridge, 1994.

\bibitem{[8]} A. Klimyk, K. Schm\"udgen:
  {\it Quantum Groups and Their Representations}.
  Springer-Verlag, Berlin Heidelberg, 1997.

\bibitem{[9]} P. \v S\v tov\'\i\v cek: Czech. J. Phys. {\bf50} (2000) 1353.

\bibitem{[10]} P. \v S\v tov\'\i\v cek: in
  {\it Proc. Internat. Sympos. Quantum Theory and Symmetries},
  Goslar 1999 (Eds. H.-D. Doebner et al.),
  World Scientific, Singapore, 2000, p.470.

\end{thebibliography}
\end{document}